\documentclass{amsart}

\usepackage{amsmath}
\usepackage{amsfonts}
\usepackage{amssymb,enumerate}
\usepackage{amsthm}
\usepackage[all]{xy}

\newtheorem{lem}{Lemma}[section]
\newtheorem{cor}[lem]{Corollary}
\newtheorem{prop}[lem]{Proposition}
\newtheorem{thm}[lem]{Theorem}

\newtheorem{Defn}[lem]{Definition}
\newtheorem{Ex}[lem]{Example}
\newtheorem{Question}[lem]{Question}
\newtheorem{Property}[lem]{Property}
\newtheorem{Properties}[lem]{Properties}
\newtheorem{Discussion}[lem]{Remark}
\newtheorem{Construction}[lem]{Construction}
\newtheorem{Subprops}{}[lem]
\newtheorem{Para}[lem]{}

\newenvironment{ex}{\begin{Ex}\rm}{\end{Ex}}

\newenvironment{subprops}{\begin{Subprops}\rm}{\end{Subprops}}
\newenvironment{para}{\begin{Para}\rm}{\end{Para}}
\newenvironment{disc}{\begin{Discussion}\rm}{\end{Discussion}}
\newenvironment{construction}{\begin{Construction}\rm}{\end{Construction}}

\newcommand{\cbc}[2]{#1(#2)}

\newcommand{\ideal}[1]{\mathfrak{#1}}
\newcommand{\m}{\ideal{m}}
\newcommand{\n}{\ideal{n}}
\newcommand{\p}{\ideal{p}}
\newcommand{\q}{\ideal{q}}
 
\newcommand{\pd}{\mathrm{pd}}

\newcommand{\depth}{\mathrm{depth}}     
\newcommand{\rank}{\mathrm{rank}}       
\newcommand{\zz}{\mathbb{Z}}

\newcommand{\injdim}{\mathrm{id}}       
\newcommand{\rhom}{\mathbf{R}\mathrm{Hom}}      
\newcommand{\lotimes}{\otimes^{\mathbf{L}}}
\newcommand{\amp}{\mathrm{amp}}
\newcommand{\HH}{\mathrm{H}}
\newcommand{\Hom}{\mathrm{Hom}} 
\newcommand{\edim}{\mathrm{edim}}

\newcommand{\fd}{\mathrm{fd}}

\newcommand{\spec}{\mathrm{Spec}}

\newcommand{\s}{\mathfrak{S}}

\newcommand{\vf}{\varphi}

\newcommand{\D}{\mathcal{D}}

\newcommand{\ol}{\overline}

\newcommand{\len}{\operatorname{length}}

\newcommand{\dist}{\operatorname{dist}}
\newcommand{\prox}{\operatorname{\sigma}}

\newcommand{\from}{\leftarrow}
\newcommand{\xra}{\xrightarrow}
\newcommand{\Graph}{\mathsf{\Gamma}}
\newcommand{\shift}{\mathsf{\Sigma}}

\newcommand{\route}{\gamma}
\newcommand{\sdc}[1]{[#1]}

\newcommand{\curv}{\operatorname{curv}}
\newcommand{\icurv}{\operatorname{inj\,curv}}
\newcommand{\da}[1]{#1^{\dagger}}
\newcommand{\tri}{\trianglelefteq}
\newcommand{\trineq}{\triangleleft}
\newcommand{\Cl}{\operatorname{Cl}}

\begin{document}

\bibliographystyle{amsplain}

\author{Anders Frankild}
\thanks{A.F., University of Copenhagen, Institute for Mathematical 
Sciences, Department of Mathematics, 
Universitetsparken 5, 2100 K\o benhavn, Denmark; email: \texttt{frankild@math.ku.dk}}

\author{Sean Sather-Wagstaff}
\thanks{S.S.-W., Department of Mathematics, California State University, Dominguez Hills, 
1000 E.~Victoria St., Carson, CA 90747 USA; email: \texttt{ssather@csudh.edu}}

\thanks{Corresponding author:  Sean Sather-Wagstaff, 
Department of Mathematics, California State University, Dominguez Hills, 
1000 E.~Victoria St., Carson, CA 90747 USA; email: \texttt{ssather@csudh.edu}}

\thanks{This research 
was conducted while 
A.F.~was funded by the Lundbeck Foundation and by
Augustinus Fonden, and S.S.-W.~was an NSF Mathematical Sciences 
Postdoctoral Research Fellow and a visitor at the University of 
Nebraska-Lincoln.}

\title[Semidualizing complexes]{The set of semidualizing 
complexes is a nontrivial metric space}

\dedicatory{Dedicated to Lars Kjeldsen, dr.~med.}

\begin{abstract}
We show that the set $\s(R)$ of 
shift-isomorphism classes of semidualizing complexes over a local 
ring $R$ admits a nontrivial metric.  We investigate the interplay between 
the metric and several algebraic operations.  
Motivated by 
the dagger duality isometry, we prove the following:
If $K,L$ are homologically bounded below and degreewise finite $R$-complexes
such that $K\lotimes_R K\lotimes_R L$ is semidualizing, then $K$ is shift-isomorphic
to $R$.
In investigating the existence of nontrivial open balls in $\s(R)$,
we prove that $\s(R)$
contains elements that are not comparable
in the reflexivity ordering if and only if it contains at least three 
distinct elements.
\\[1mm]
\end{abstract}
\maketitle

\section*{Introduction} 

Much research in commutative algebra
is devoted to duality.
One example of this is the work of Grothendieck and
Hartshorne~\cite{hartshorne:rad} which includes an investigation of 
the duality properties of 
finite modules and complexes with respect to a dualizing complex.
A second example is the work of Auslander and 
Bridger~\cite{auslander:adgeteac,auslander:smt} where a 
class of modules is identified, those of finite G-dimension,
having good duality properties with respect to the 
ring.

These examples are antipodal in the sense that 
each one is devoted to the reflexivity 
properties of finite modules and complexes with respect to a 
semidualizing complex.  
See~\ref{homothety} for precise definitions.
Examples of semidualizing complexes include 
the ring itself and the dualizing complex, if it exists.
Another useful example is
the dualizing complex of a 
local homomorphism of finite G-dimension,
as constructed by Avramov and Foxby~\cite{avramov:rhafgd}.  
The study of the general situation was initiated by 
Foxby~\cite{foxby:gmarm},
Golod~\cite{golod:gdagpi}, and
Christensen~\cite{christensen:scatac}.
 
We denote by $\s(R)$ the set of shift-isomorphism classes 
of semidualizing 
complexes over a local ring $R$;  the  class
of a given semidualizing complex $K$ is denoted $\sdc{K}$.
The work in the current paper is part of an ongoing research effort 
on our part to analyze
the structure of the set $\s(R)$ in its entirety.
That $\s(R)$ has more structure than other collections of 
complexes is demonstrated by the fact that one can inflict 
upon $\s(R)$ an ordering given by reflexivity;
see~\ref{homothety} and~\ref{para100}.  
Further structure is demonstrated in~\cite{sather:divisor} where it is 
observed that, when $R$ is a Cohen-Macaulay normal domain,
the set 
$\s(R)$ is naturally a 
subset of the divisor class group $\Cl(R)$.  The analysis of this inclusion 
yields, for instance, a complete description of $\s(R)$ for certain 
classes of rings; see, e.g., Example~\ref{nontrivial}.

The main idea in the present work is to use numerical data from 
the complexes in $\s(R)$ that are comparable under the ordering to 
give a measure of their proximity.  The distance between two arbitrary 
elements $\sdc{K},\sdc{L}$ of $\s(R)$ is then described via chains of 
pairwise comparable elements starting with $\sdc{K}$ and ending with 
$\sdc{L}$.  Details of the 
construction and its basic properties are given in 
Section~\ref{sec2}.  One main result, advertised in the title, 
is  contained in Theorem~\ref{prop1}.

\medskip

\noindent \textbf{Theorem A.}  
\textit{The set $\s(R)$ is a metric space.}

\medskip

Theorem~\ref{prop201}, stated next, 
shows that the metric is not 
equivalent to the trivial one,
unless $\s(R)$ is itself almost trivial.
It also implies that $\s(R)$ quite frequently contains elements
that are
noncomparable in the ordering.

\medskip

\noindent \textbf{Theorem B.}  
\textit{For a local ring $R$ the following conditions are equivalent:
\begin{enumerate}[\quad\rm(i)]
\item \label{item207'}
There exist elements of $\s(R)$ that are not comparable.
\item \label{item208'}
$\s(R)$ has cardinality at least 3.
\item \label{item209'}
There exists $\sdc{K}\in\s(R)$ and $\delta>0$ such that the open 
ball $B(\sdc{K},\delta)$ satisfies $\{\sdc{K}\}\subsetneq 
B(\sdc{K},\delta)\subsetneq\s(R)$.
\end{enumerate}}

\medskip

This result follows in part from an analysis motivated by 
Proposition~\ref{prop4}:  If $R$ admits a dualizing complex $D$, then the map
$\s(R)\to\s(R)$ given by
$\sdc{K}\mapsto\sdc{\rhom_R(K,D)}$ is an isometric involution.
This fact led us to investigate the fixed points 
of this involution.  Corollary~\ref{fixed02}
shows that the existence of such a fixed point implies that $R$ is Gorenstein;
it is a consequence of Theorem~\ref{fixed01}, stated next.

\medskip

\noindent \textbf{Theorem C.}  
\textit{Let $R$ be a local ring and $K,L$ homologically 
bounded below and degreewise finite 
$R$-complexes. 
If $K\lotimes_RK\lotimes_R L$ is 
semidualizing, then $K$ is shift-isomorphic 
to $R$ in the derived category $\D(R)$.
}

\medskip

Section~\ref{sec3} describes the behavior of the metric with respect 
to change of rings along a local homomorphism of finite flat dimension.
In particular, these operations give a recipe for constructing noncomparable 
semidualizing complexes;  see Corollary~\ref{noncom}.
We conclude with Section~\ref{sec4}, which
consists of explicit computations.

\section{Complexes} \label{sec1}

This section consists of background and includes most of the 
definitions and notational conventions used throughout the rest of 
this paper.

\medskip
\noindent\emph{Throughout this work, $(R,\m,k)$ and $(S,\n,l)$ are local Noetherian 
commutative rings.}

\medskip
An \emph{$R$-complex} is a sequence of $R$-module homomorphisms
\[ X=\cdots\xra{\partial^X_{i+1}}X_i\xra{\partial^X_{i}}X_{i-1} 
\xra{\partial^X_{i-1}}\cdots \]
with $\partial_i^X\partial_{i+1}^X=0$ for each $i$.  We work 
in the derived category $\D(R)$ whose objects are the $R$-complexes;   
references on the subject
include~\cite{gelfand:moha,hartshorne:rad,verdier:cd,verdier:1}.
For $R$-complexes $X$ and $Y$ 
the left derived tensor product complex 
is denoted $X\lotimes_R Y$ and 
the right derived homomorphism complex is 
$\rhom_R(X,Y)$.  For an integer $n$, 
the $n$th \emph{shift} or \emph{suspension} of $X$ is denoted
$\shift^n X$ where $(\shift^n 
X)_i=X_{i-n}$ and $\partial_i^{\shift^n X}=(-1)^n\partial_{i-n}^X$.
The symbol ``$\simeq$'' indicates an
isomorphism in $\D(R)$ and ``$\sim$'' indicates an isomorphism up to shift.

The \emph{infimum}, \emph{supremum}, and \emph{amplitude} of a complex 
$X$ are, respectively, 
\begin{align*}
\inf(X)&=\inf\{i\in\zz\mid\HH_i(X)\neq 0\} \\
\sup(X)&=\sup\{i\in\zz\mid\HH_i(X)\neq 0\} \\
\amp(X)&=\sup(X)-\inf(X) 
\end{align*}
with the conventions $\inf\emptyset=\infty$ 
and $\sup\emptyset=-\infty$.
The complex $X$ is \emph{homologically finite}, respectively 
\emph{homologically degreewise finite}, if its total 
homology module $\HH(X)$, respectively each individual homology 
module $\HH_i(X)$, is a finite $R$-module.

The \emph{$i$th Betti number} and \emph{Bass number} of 
a homologically finite complex of $R$-modules $X$ are, respectively,
\[ \beta_i^R(X)=\rank_k(\HH_{-i}(\rhom_R(X,k)))
\quad\text{and}\quad
\mu^i_R(X)=\rank_k(\HH_{-i}(\rhom_R(k,X)). \]
The \emph{Poincar\'{e} series} and \emph{Bass series}
of $X$ are the formal Laurent series
\[ P^R_X(t)=\sum_{i\in\zz}\beta_i^R(X) t^{i}
\qquad \text{and} \qquad
I_R^X(t)=\sum_{i\in\zz}\mu^i_R(X) t^i. \]
The projective, injective, and flat dimensions of $X$ are denoted 
$\pd_R(X)$, $\injdim_R(X)$, and $\fd_R(X)$, respectively;  
see~\cite{foxby:hacr}.

The Bass series of a local homomorphism of finite flat dimension 
is an important invariant that will appear in several contexts in this 
work.

\begin{para}  \label{paraBass}
A ring homomorphism $\vf\colon R\to S$ is local when 
$\vf(\m)\subseteq\n$.  In this event, the flat dimension of 
$\vf$ is defined as $\fd(\vf)=\fd_R(S)$, and the 
\emph{depth} of $\vf$ is $\depth(\vf)=\depth(S)-\depth(R)$.
When $\fd(\vf)$ is finite,
the \emph{Bass series} of $\vf$ is the formal Laurent 
series with nonnegative integer coefficients $I_{\vf}(t)$ satisfying 
the formal equality
\[ I_{S}^{S}(t)=I_R^R(t)I_{\vf}(t).\]
The existence of $I_{\vf}(t)$
is given by~\cite[(5.1)]{avramov:bsolrhoffd} 
or~\cite[(7.1)]{avramov:rhafgd}.
The 
map $\vf$ is \emph{Gorenstein} at $\n$ if $I_{\vf}(t)=t^d$ for some 
integer $d$ (in which case $d=\depth(\vf)$) equivalently, if $I_{\vf}(t)$ is a Laurent
polynomial.
\end{para}

Our metric utilizes the curvature of 
a homologically finite complex, as introduced by 
Avramov~\cite{avramov:mwer}.  It provides an exponential measure of 
the growth of the Betti numbers of the complex.
Let $F(t)=\sum_{n\in\zz}a_nt^n$ be
a formal Laurent series with nonnegative integer coefficients.
The \emph{curvature} 
of $F(t)$ is
\[ \curv(F(t))=\limsup_{n\to\infty}\sqrt[n]{a_n}. \]
Of the following properties, parts~\eqref{itemcurv} 
and~\eqref{item01} follow  from the definition.  For 
part~\eqref{item02}, argue as in the proof 
of~\cite[(4.2.4.6)]{avramov:ifr}.

\begin{para} \label{product}
Let $F(t),G(t)$ be formal Laurent series with nonnegative integer 
coefficients.
\begin{enumerate}[\quad(a)]
\item \label{itemcurv}
For each integer $d$, there is an equality $\curv(F(t))=\curv(t^dF(t))$.
\item  \label{item01}
A coefficientwise inequality $F(t)\preccurlyeq 
G(t)$ implies 
$\curv(F(t))\leq\curv(G(t))$.
\item  \label{item02} There is an equality
$\curv(F(t)G(t))=\max\{\curv(F(t)),\curv(G(t))\}$.
\end{enumerate}
\end{para}

\begin{para} \label{paraRadii}
Let $X$ be a homologically degreewise finite $R$-complex.
The \emph{curvature} and \emph{injective curvature} of $X$ are
\[ \curv_R(X)=\curv(P^R_X(t)) \qquad \text{and} \qquad
\icurv_R(X)=\curv(I_R^X(t)).\]
For a local homomorphism of finite 
flat dimension $\vf$, the \emph{injective curvature} 
of $\vf$ is 
\[ \icurv(\vf)=\curv(I_{\vf}(t)).\]
In particular, the map $\vf$ is Gorenstein at $\n$ if and only if $\icurv(\vf)=0$.
\end{para}

\begin{para}  \label{para8}
Let $X,Y$ be homologically finite complexes of $R$-modules.
\begin{enumerate}[\quad\rm(a)]
\item \label{item10}
If $\vf\colon R\to S$ is a local homomorphism, then 
$\curv_{S}(X\lotimes_R S)=\curv_R(X)$.
\item \label{item12}
There are inequalities
$0\leq\curv_R(X)<\infty$.
\item \label{item13}
The following 
conditions are equivalent:
\begin{enumerate}[\quad\rm(i)]
\item $\curv_R(X)<1$,
\item $\curv_R(X)=0$,
\item $\pd_R(X)$ is finite.
\end{enumerate}
\end{enumerate}
\end{para}

\begin{proof}
The formal equality
$P^S_{X\lotimes_R S}(t)=P^R_X(t)$ is 
by~\cite[(1.5.3)]{avramov:rhafgd}, and part~\eqref{item10} follows. 
Apply~\cite[(4.1.9),(4.2.3.5),(4.2.3.1)]{avramov:ifr} to a truncation of 
the minimal free resolution of $X$
to verify~\eqref{item12} and~\eqref{item13}. 
\end{proof}

Next, we turn to
semidualizing complexes and their reflexive objects.

\begin{para} \label{homothety}
For homologically finite $R$-complexes $K$ and $X$ one has 
natural homothety and biduality homomorphisms, respectively.  
\begin{align*}
\chi^R_K & \colon R\to\rhom_R(K,K) \\
\delta^K_X & \colon X\to\rhom_R(\rhom_R(X,K),K)
\end{align*}
The complex $K$ is \emph{semidualizing} if $\chi^R_K$ is an isomorphism; e.g.,
$R$ is 
semidualizing.  

A complex $D$ is \emph{dualizing} if it is 
semidualizing and has finite injective dimension;  
see~\cite[Chap.~V]{hartshorne:rad}. 
Dualizing complexes are unique up to shift-isomorphism.
Any
homomorphic image of a local Gorenstein ring, e.g., any complete 
local ring, 
admits a dualizing complex by~\cite[(V.10.4)]{hartshorne:rad}.  
When $D$ is dualizing for $R$, one has
$I_R^D(t)=t^d$ for some integer $d$ by~\cite[(V.3.4)]{hartshorne:rad}.

When $K$ is semidualizing, the complex $X$ is \emph{$K$-reflexive} 
if  $\rhom_R(X,K)$ is homologically bounded and  
$\delta^K_X$ is an isomorphism; e.g., $R$ 
and $K$ are $K$-reflexive.
When $R$ admits a dualizing complex $D$, each homologically finite 
complex $X$ is $D$-reflexive by~\cite[(V.2.1)]{hartshorne:rad}.  A 
complex is $R$-reflexive exactly when it has finite 
G-dimension by~\cite[(2.3.8)]{christensen:gd}. 
\end{para}

The Poincar\'{e} and Bass series of a semidualizing complex 
are linked by~\cite[(1.5.3)]{avramov:rhafgd}.

\begin{para} \label{para101}
When $K$ is a semidualizing $R$-complex, 
there is a formal equality
\[ P^R_K(t)I_R^K(t)=I^R_R(t). \]
\end{para}

Here is the fundamental object of study in this work.

\begin{para} \label{para100}
The set 
of shift-isomorphism classes of semidualizing $R$-complexes is 
denoted $\s(R)$.  
The class in $\s(R)$ of a semidualizing 
complex $K$  is denoted $\sdc{K}$.
For $\sdc{K},\sdc{L}\in \s(R)$ write $\sdc{K}\tri \sdc{L}$ if 
$L$ is $K$-reflexive; this is independent of the representatives for $\sdc{K}$ 
and $\sdc{L}$, and
$\sdc{K}\tri \sdc{R}$.  If $D$ is dualizing, then
$\sdc{D}\tri\sdc{L}$.
\end{para}

\begin{para} \label{para0}
If $\sdc{K},\sdc{L}\in\s(R)$ and $\sdc{K}\tri\sdc{L}$, then
$\rhom_R(L,K)$ is 
semidualizing and $K$-reflexive,
that is,
$\sdc{K}\tri \sdc{\rhom_R(L,K)}$;  
see~\cite[(2.11)]{christensen:scatac}.  In particular, if $D$ is dualizing for $R$, then
the complex 
$\da{L}=\rhom_R(L,D)$
is semidualizing; there are equalities
$I^{L^{\dagger}}_R(t)=t^{d}P^R_L(t)$
and 
$P_{L^{\dagger}}^R(t)=t^{-d}I_R^L(t)$
for some $d\in\mathbb{Z}$ by~\cite[(1.7.7)]{christensen:scatac},
and so $\icurv_R(L)=\curv_R(L^{\dagger})$ by~\ref{product}\eqref{itemcurv}.
\end{para}

\begin{para} \label{gerko1}
For semidualizing complexes $K,L,M$, 
consider the composition morphism
\[ \xi_{MLK}\colon\rhom_R(M,L)\lotimes_R\rhom_R(L,K)\to\rhom_R(M,K).\]
This is an isomorphism when
$L$ and $M$ are $K$-reflexive
and $M$ is $L$-reflexive by~\cite[(3.3)]{gerko:sdc}, 
and a   
formal equality of Laurent series follows 
from~\cite[(1.5.3)]{avramov:rhafgd}
\[ P^R_{\rhom_R(M,L)}(t)P^R_{\rhom_R(L,K)}(t)=P^R_{\rhom_R(M,K)}(t). \]
In particular, when $M=R$ 
the morphism is of the form
$L\lotimes_R \rhom_R(L,K)\to K$,
and when $L$ is $K$-reflexive 
one has $P^R_L(t)P^R_{\rhom_R(L,K)}(t)=P^R_K(t)$.
\end{para}

\begin{para} \label{paraRadius}
For $\sdc{K}$ in $\s(R)$, the 
quantities $\curv_R(K)$ and $\icurv_R(K)$
are well-defined.
\label{para6}
There are inequalities
\[ 0\leq\curv_R(K)\leq\icurv_R(R)<\infty
\quad \text{and} \quad
0\leq\icurv_R(K)\leq\icurv_R(R)<\infty\]
and  the following 
conditions are equivalent:
\begin{enumerate}[\quad\rm(i)]
\item $\curv_R(K)<1$,
\item $\curv_R(K)=0$,
\item $\sdc{K}= \sdc{R}$. 
\end{enumerate}
\end{para}

\begin{proof}
For the first statement see~\ref{product}\eqref{itemcurv}, 
while the equivalence of (i)--(iii) follows from~\ref{para8}\eqref{item13} 
and~\cite[(8.1)]{christensen:scatac}.
For the inequalities,
pass to the completion of $R$ to assume that $R$ admits a dualizing 
complex $D$.  With $(-)^{\dagger}$ as in~\ref{para0}, 
use~\ref{homothety} and~\ref{gerko1} to 
verify the equality
$P^R_D(t)=P^R_K(t)P^R_{\da{K}}(t)$. 
 With~\ref{product}\eqref{item02}
 this provides the first equality in the next sequence
 while the second is in~\ref{para0}.
 \[ \curv_R(K)\leq
\max\{\curv_R(K),\curv_R(\da{K})\} =\curv_R(D) 
=\icurv_R(R) \]
This gives one of the inequalities, while the others follow from~\ref{para8}\eqref{item12}
and~\ref{para0}.  
\end{proof}

\section{The metric} \label{sec2}

Here is the first step of the construction of the metric on $\s(R)$.

\begin{para} \label{para7}
For $\sdc{K},\sdc{L}$ in $\s(R)$ with $\sdc{K}\tri \sdc{L}$, set 
\[ \prox_R(\sdc{K},\sdc{L})=\curv_R(\rhom_R(L,K)).\]
\end{para}

Apply~\ref{gerko1} and~\ref{para6} to $\sdc{\rhom_R(L,K)}$ in order to establish the following.

\begin{lem}  \label{para1}
For $\sdc{K},\sdc{L}\in\s(R)$ with $\sdc{K}\tri 
\sdc{L}$, the quantity 
$\prox_R(\sdc{K},\sdc{L})$ is well-defined and nonnegative.
Furthermore, the following conditions are equivalent:
\begin{enumerate}[\quad\rm(i)]
\item $\prox_R(\sdc{K},\sdc{L})<1$,
\item $\prox_R(\sdc{K},\sdc{L})=0$,
\item $\sdc{K} =\sdc{L}$.\qed
\end{enumerate}
\end{lem}

The following simple
construction helps us visualize the metric. 

\begin{construction}
Let $\Graph(R)$ be the 
directed graph whose 
vertex set is $\s(R)$ and whose directed edges $\sdc{K}\to \sdc{L}$ correspond 
exactly to the inequalities $\sdc{K}\tri \sdc{L}$.  
Graphically, ``smaller'' semidualizing 
modules will be drawn below ``larger'' ones. 
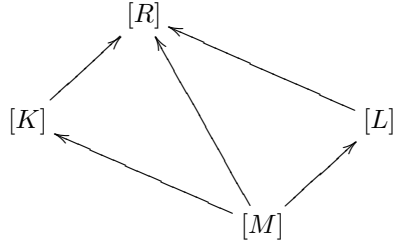
\[ \xymatrix{
 & \sdc{R} & \\
\sdc{K} \ar[ur] & & & \sdc{L} \ar[ull] \\
& & \sdc{M} \ar[ull] \ar[ur] \ar[uul]
} \]
\end{construction}

The metric will arise from the graph 
$\Graph(R)$ with a ``taxi-cab metric'' in mind where $\prox_R$ is used 
to measure the length of the edges.

\begin{para} \label{para3}
A \emph{route} $\route$ from $\sdc{K}$ to $\sdc{L}$ in $\Graph(R)$ is 
a subgraph of $\Graph(R)$ of the form 
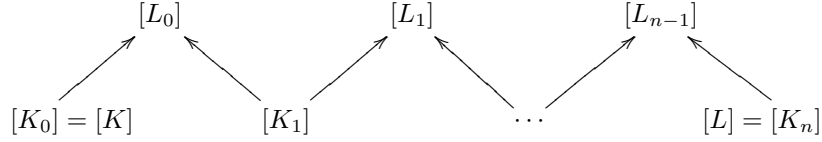
\[ \xymatrix{
  & \sdc{L_0} & & \sdc{L_1} & & \sdc{L_{n-1}} \\ 
   \sdc{K_0}=\sdc{K} \hspace{-1cm} 
 \ar[ur] & & \sdc{K_1} \ar[ul]\ar[ur] & & 
 \cdots \ar[ul]\ar[ur] & & \hspace{-1cm} \sdc{L}=\sdc{K_n}    
 \ar[ul] 
 } \]
and the \emph{length} of the route $\route$ is the sum of the lengths
of its edges
\[ \len_R(\route) 
=\prox_R(\sdc{K_0},\sdc{L_0})+
\prox_R(\sdc{K_1},\sdc{L_0})+\cdots+
\prox_R(\sdc{K_n},\sdc{L_{n-1}}). \]
By Lemma~\ref{para1}, there is an inequality $\len_R(\route)\geq 
0$.
\end{para}

\begin{disc} \label{routes}
The fact that $\Graph(R)$ is a \emph{directed} graph is only used to 
keep track of routes in $\Graph(R)$.
We define the metric in terms of routes 
instead of arbitrary paths in order to keep the 
notation simple.  For instance, the proof of Theorem~\ref{prop2} 
would be even more notationally complicated without the directed 
structure.
Note that
the metric that arises by considering arbitrary 
paths in $\Graph(R)$ is equal to the one we construct below.  Indeed,
any path in $\Graph(R)$ from $\sdc{K}$ to $\sdc{L}$ can be expressed as a route
of the same length by 
inserting trivial edges $\sdc{M}\to \sdc{M}$.
\end{disc}

Here are some specific routes whose lengths will give rise to 
bounds on the metric.

\begin{para} \label{para4}
Since $\sdc{K}\tri \sdc{R}$ and $\sdc{L}\tri \sdc{R}$, 
a route $\route_1$ from $\sdc{K}$ to $\sdc{L}$ always exists
\[ \xymatrix{
 & \sdc{R} \\
\sdc{K} \ar[ur] & & \sdc{L} \ar[ul]
} \]
with
$\len_R(\route_1)=\curv_R(K)+\curv_R(L)$.
In particular, the graph $\Graph(R)$ is connected.  We shall see 
in Theorem~\ref{prop201} 
that the graph is not complete in general.

When $R$ admits a dualizing complex $D$, another route $\route_2$ from 
$\sdc{K}$ to $\sdc{L}$ 
is 
\[ \xymatrix{
 & \sdc{K} & &  \sdc{L} \\
\sdc{K} \ar[ur] & & \sdc{D} \ar[ur]\ar[ul] & & \sdc{L} \ar[ul]
} \]
and
$\len_R(\route_2)=\curv_R(K^{\dagger})+\curv_R(L^{\dagger})
=\icurv_R(K)+\icurv_R(L)$ by~\ref{para0}.
\end{para}

The next properties are straightforward to verify.

\begin{para}
Fix $\sdc{K},\sdc{L},\sdc{M}$ in $\s(R)$.

\begin{subprops} \label{subprop1}
The set of routes from $\sdc{K}$ to $\sdc{L}$ is 
in length-preserving bijection with the set of routes from $\sdc{L}$ to 
$\sdc{K}$.  
\end{subprops}

\begin{subprops} \label{subprop2}
The diagram $\sdc{K}\to \sdc{K}\from \sdc{K}$ gives a route from 
$\sdc{K}$ to $\sdc{K}$ with length 0.  
\end{subprops}

\begin{subprops} \label{subprop3}
Let $\route$ be a route from $\sdc{K}$ to $\sdc{L}$ and $\route'$ a route 
from $\sdc{L}$ to $\sdc{M}$.
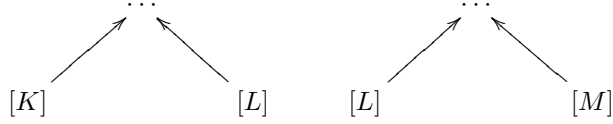
\[ \xymatrix{
 & \cdots & & & \cdots \\
\sdc{K} \ar[ur] & & \sdc{L}\ar[ul] & \sdc{L}\ar[ur] & & \sdc{M}\ar[ul] 
} \]
Let $\route\route'$ denote the concatenation of $\route$ and $\route'$
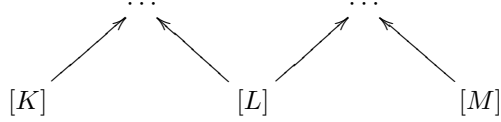
\[ \xymatrix{
 & \cdots & & \cdots \\
\sdc{K} \ar[ur] & & \sdc{L}\ar[ur]\ar[ul] & & \sdc{M}\ar[ul] 
} \]
It is immediate that
$\len_R(\route\route')=\len_R(\route)+\len_R(\route')$.
\end{subprops}
\end{para}

Here is the definition of our metric on $\s(R)$.

\begin{para} \label{defn}
The \emph{distance} from $\sdc{K}$ to $\sdc{L}$ in $\s(R)$ is
\[ \dist_R(\sdc{K},\sdc{L})=\inf\{\len_R(\route)\mid\text{$\route$ is a route 
from $\sdc{K}$ to $\sdc{L}$ in $\Graph(R)$}\}.\]
\end{para}

The next result  is Theorem A from the introduction.

\begin{thm} \label{prop1}
The function $\dist_R$ is a metric on $\s(R)$.
\end{thm}

\begin{proof}
Fix $\sdc{K},\sdc{L}$ in $\s(R)$.
The inequality $\dist_R(\sdc{K},\sdc{L})\geq 0$ 
is satisfied since $\len_R(\route)\geq 0$ for each
route $\route$ from $\sdc{K}$ to $\sdc{L}$ in $\Graph(R)$
by~\ref{para3}, and
at least one such route 
exists by~\ref{para4}.  With this, the computation 
in~\ref{subprop2} 
shows that 
$\dist_R(\sdc{K},\sdc{K})=0$.  
If $\dist_R(\sdc{K},\sdc{L})=0$, 
then there is 
a route $\route$ from $\sdc{K}$ to $\sdc{L}$ 
in $\Graph(R)$ with  $\len_R(\route)<1$.  
Using the 
notation for $\route$ as in~\ref{para3}, one has 
\[ \prox_R(\sdc{K_i},\sdc{L_j})<1 \qquad
\text{for $j=0,\ldots,n-1$ and $i=j,j+1$} \]
and therefore by Lemma~\ref{para1} there are equalities
$\sdc{K_i}=\sdc{L_j}$ and so $\sdc{K}=\sdc{L}$.
Thus, $\dist_R(\sdc{K},\sdc{L})\geq 0$ with equality if and only if
$\sdc{K}= \sdc{L}$.

It follows from~\ref{subprop1} that 
$\dist_R(\sdc{K},\sdc{L})=\dist_R(\sdc{L},\sdc{K})$.  
To verify the triangle inequality, fix $\sdc{M}$ in 
$\s(R)$.  For each real number $\epsilon>0$, we
will verify the inequality
\begin{equation} \label{triangle} \tag{$\dagger$}
\dist_R(\sdc{K},\sdc{M})<
\dist_R(\sdc{K},\sdc{L})+\dist_R(\sdc{L},\sdc{M}) +\epsilon
\end{equation}
and the inequality
$\dist_R(\sdc{K},\sdc{M})\leq
\dist_R(\sdc{K},\sdc{L})+\dist_R(\sdc{L},\sdc{M})$
will follow.  Fix an $\epsilon>0$ and choose routes
$\route$ from $\sdc{K}$ to $\sdc{L}$ 
and $\route'$ from $\sdc{L}$ to $\sdc{M}$ 
with 
\[ \len_R(\route)<\dist_R(\sdc{K},\sdc{L})+\epsilon/2
\qquad\text{and}\qquad 
\len_R(\route')<\dist_R(\sdc{L},\sdc{M})+\epsilon/2; \]  
such routes exist by the basic properties of the infimum.
The concatenation 
$\route\route'$ is 
a route from $\sdc{K}$ to $\sdc{M}$, explaining (1) in 
the following sequence 
\begin{align*}
\dist_R(\sdc{K},\sdc{M})
&\stackrel{(1)}{\leq}\len_R(\route\route') \\
&\stackrel{(2)}{=}\len_R(\route)+\len_R(\route') \\
&\stackrel{(3)}{<}\dist_R(\sdc{K},\sdc{L})+\dist_R(\sdc{L},\sdc{M})+\epsilon
\end{align*}
while (2)
is by~\ref{subprop3}, and (3) is by the choice of 
$\route$ and $\route'$.
\end{proof}

\begin{disc}
Given $\sdc{K},\sdc{L}$ in $\s(R)$, it is not clear
from the definition of the metric that there exists a route $\route$ 
from $\sdc{K}$ to $\sdc{L}$ such that 
$\dist_R(\sdc{K},\sdc{L})=\len_R(\route)$.  If $\s(R)$ 
is finite, more generally, if the set 
$\{\curv_R(M)\mid\sdc{M}\in\s(R)\}$ 
is finite, such a route 
would exist.  (Compare this to~\cite[Problem 4.3.8]{avramov:ifr} which
asks whether the curvature function takes on finitely many values in 
total.)  The next result gives one 
criterion guaranteeing that such a route exists:  when 
$\sdc{K}\tri\sdc{L}$, the trivial route 
$\sdc{K}\to\sdc{L}\from\sdc{L}$ is length-minimizing.
\end{disc}

\begin{thm} \label{prop2}
For $\sdc{K}\tri \sdc{L}$ in $\s(R)$, one has 
$\dist_R(\sdc{K},\sdc{L})=\prox_R(\sdc{K},\sdc{L})$.
\end{thm}

\begin{proof}
The route $\sdc{K}\to \sdc{L}\from \sdc{L}$ has length 
$\prox_R(\sdc{K},\sdc{L})$
giving the inequality
$\dist_R(\sdc{K},\sdc{L})\leq\prox_R(\sdc{K},\sdc{L})$.
Fix a route $\route$ 
from $\sdc{K}$ to $\sdc{L}$ in $\Graph(R)$. 
We verify the inequality 
$\prox_R(\sdc{K},\sdc{L})\leq\len_R(\route)$;  
this will yield the inequality
$\prox_R(\sdc{K},\sdc{L})\leq\dist_R(\sdc{K},\sdc{L})$, completing 
the proof.
With notation for $\route$ as 
in~\ref{para3}, set
\[ 
P_{i,j}(t)=P^R_{\rhom_R(L_j,K_i)}(t)=
P^R_{K_i}(t)/P^R_{L_j}(t)
\quad \text{for $j=0,\ldots,n-1$ and $i=j,j+1$} \]
where the second equality is from~\ref{gerko1}.
This gives (1) and (6) in the following sequence where the
formal equalities 
hold in the field of fractions of the ring of formal Laurent 
series with integer coefficients.
\begin{align*}
P_{0,0}(t)P_{1,0}(t)P_{1,1}(t)\cdots P_{n,n-1}(t) \hspace{-2cm} \\
& \stackrel{(1)}{=} 
  \frac{P^R_{K_0}(t)}{P^R_{L_0}(t)} 
  \frac{P^R_{K_1}(t)}{P^R_{L_0}(t)}
  \frac{P^R_{K_1}(t)}{P^R_{L_1}(t)} \cdots 
  \frac{P^R_{K_{n-1}}(t)}{P^R_{L_{n-1}}(t)} 
  \frac{P^R_{K_n}(t)}{P^R_{L_{n-1}}(t)} \\
& \stackrel{(2)}{=} 
  \frac{P^R_{K_0}(t)}{P^R_{L_0}(t)} 
  \frac{P^R_{K_1}(t)}{P^R_{L_0}(t)}
  \frac{P^R_{K_1}(t)}{P^R_{L_1}(t)} \cdots 
  \frac{P^R_{K_{n-1}}(t)}{P^R_{L_{n-1}}(t)} 
  \frac{P^R_{K_n}(t)}{P^R_{L_{n-1}}(t)} 
  \frac{P^R_{K_n}(t)}{P^R_{K_n}(t)} \\
& \stackrel{(3)}{=} 
  P^R_{K_0}(t)
  \biggl[\frac{P^R_{K_1}(t)}{P^R_{L_0}(t)}\biggr]^2
  \biggl[\frac{P^R_{K_2}(t)}{P^R_{L_1}(t)}\biggr]^2
  \cdots
  \biggl[\frac{P^R_{K_n}(t)}{P^R_{L_{n-1}}(t)}\biggr]^2
  \frac{1}{P^R_{K_n}(t)} \\
& \stackrel{(4)}{=} 
  \biggl[\frac{P^R_{K_0}(t)}{P^R_{K_n}(t)}\biggr]
  \biggl[\frac{P^R_{K_1}(t)}{P^R_{L_0}(t)}\biggr]^2
  \biggl[\frac{P^R_{K_2}(t)}{P^R_{L_1}(t)}\biggr]^2
  \cdots
  \biggl[\frac{P^R_{K_n}(t)}{P^R_{L_{n-1}}(t)}\biggr]^2 \\
& \stackrel{(5)}{=} 
  \biggl[\frac{P^R_{K}(t)}{P^R_{L}(t)}\biggr]
  \biggl[\frac{P^R_{K_1}(t)}{P^R_{L_0}(t)}\biggr]^2
  \biggl[\frac{P^R_{K_2}(t)}{P^R_{L_1}(t)}\biggr]^2
  \cdots
  \biggl[\frac{P^R_{K_n}(t)}{P^R_{L_{n-1}}(t)}\biggr]^2 \\
& \stackrel{(6)}{=} 
   P^R_{\rhom_R(L,K)}(t)
   \prod_{i=1}^n P_{i,i-1}(t)^2 \\
& \stackrel{(7)}{\succcurlyeq} P^R_{\rhom_R(L,K)}(t) t^d
\end{align*}
Here $d$ is twice the sum 
of the orders of the Laurent 
series $P_{i,i-1}(t)$.
Equality (2) is trivial, (3) and (4) are obtained
by rearranging the factors, (5) is by the choice of $K_0$ 
and $K_n$, and (7) follows from the 
fact that the coefficients of each $P_{i,i-1}(t)$ 
are nonnegative integers.
With~\ref{product}\eqref{item01} this
explains (11) in the following sequence
\begin{align*}
\len_R(\route) 
& \stackrel{(8)}{=}
  \curv(P_{0,0}(t))+\curv(P_{1,0}(t)) +\curv(P_{1,1}(t)) +\cdots 
  +\curv(P_{n,n-1}(t)) \\
& \stackrel{(9)}{\geq}
  \max\{\curv(P_{0,0}(t)),\curv(P_{1,0}(t)),\curv(P_{1,1}(t)),\cdots, 
  \curv(P_{n,n-1}(t))\} \hspace{-1cm} \\
& \stackrel{(10)}{=}
  \curv(P_{0,0}(t) P_{1,0}(t)P_{1,1}(t)\cdots P_{n,n-1}(t)) \\
& \stackrel{(11)}{\geq}
  \curv(P^R_{\rhom_R(L,K)}(t)) \\
& \stackrel{(12)}{=}
  \prox_R(\sdc{K},\sdc{L})
\end{align*}
where (8) and (12) are by definition, 
(9) is by the nonnegativity of each $\curv(P_{i,j}(t))$, and
(10) is by~\ref{product}\eqref{item02}. 
This completes the proof.
\end{proof} 

The computations in~\ref{para6} and~\ref{para4} provide 
bounds on the metric.

\begin{prop} \label{para5}
For $\sdc{K}$ and $\sdc{L}$ in $\s(R)$, there are inequalities
\[ \dist_R(\sdc{K},\sdc{L})\leq\curv_R(K)+\curv_R(L)\leq 
2\icurv_R(R).\]
In particular, the metric is completely bounded.
Furthermore, there are inequalities
\[ \dist_R(\sdc{K},\sdc{L})\leq\icurv_R(K)+\icurv_R(L)
\leq 
2\icurv_R(R) \]
when $R$ admits a dualizing complex. \qed
\end{prop}

\begin{disc} \label{topology}
The topology on $\s(R)$ induced by the metric
is trivial.  Indeed, 
Lemma~\ref{para1} and Theorem~\ref{prop2} imply that 
the singleton $\{\sdc{K}\}$ is exactly
the open ball of radius 1 centered at 
the point $\sdc{K}$.  Similarly, using the 
upper bound established in Proposition~\ref{para5}, the open ball of radius 
$2\icurv_R(R)+1$ is $\s(R)$ itself.  On the other hand, in Theorem~\ref{prop201} we 
show that, if $\s(R)$ contains at least three elements, then
$\s(R)$ has nontrivial open balls.
\end{disc}

\section{Dagger duality, fixed points, and nontriviality of the metric} \label{sec6}

The first result of this section uses notation from~\ref{para0}.

\begin{prop} \label{prop4}
If $R$ admits a dualizing complex $D$, then the map
$\Delta\colon\s(R)\to\s(R)$ given by sending $\sdc{K}$ to
$\sdc{K^{\dagger}}$ is an isometric involution of $\s(R)$.
\end{prop}

\begin{proof}
The map $\Delta$ is an involution of $\s(R)$ by~\ref{homothety}.  To show that it is an 
isometry, it suffices to verify the following containment of subsets 
of $\mathbb{R}$.
\begin{equation}
\{\len_R(\route)\mid\text{$\route$ a route $\sdc{K}$ to 
$\sdc{L}$}\} \subseteq
\{\len_R(\route_1)\mid\text{$\route_1$ a route $\sdc{L^{\dagger}}$ to 
$\sdc{K^{\dagger}}$}\}
\tag{$\ddagger$} 
\label{eqn2}
\end{equation}
Indeed, this will give the inequalities in the following sequence
\[ \dist_R(\sdc{K},\sdc{L})=\dist_R(\sdc{K^{\dagger\dagger}},\sdc{L^{\dagger\dagger}})
\leq\dist_R(\sdc{K^{\dagger}},\sdc{L^{\dagger}})\leq\dist_R(\sdc{K},\sdc{L}) \]
while~\ref{homothety} explains the equality; thus, equality is forced
at each step.  

When $\sdc{K}\tri \sdc{L}$, one concludes from~\cite[(3.9)]{frankild:appx}
that
$\sdc{L^{\dagger}}\tri \sdc{K^{\dagger}}$.  Furthermore, there is a sequence 
of equalities
where the middle equality
is by the isomorphism
$\rhom_R(K^{\dagger},L^{\dagger})\simeq\rhom_R(L,K)$
in~\cite[(1.7(a))]{frankild:appx}.
\[ \prox_R(\sdc{L^{\dagger}},\sdc{K^{\dagger}})
=  \curv_R(\rhom_R(K^{\dagger},L^{\dagger}))
=  \curv_R(\rhom_R(L,K)) 
=  \prox_R(\sdc{K},\sdc{L}) \]

To verify~\eqref{eqn2}, 
let $\route$ be a route from $\sdc{K}$ to $\sdc{L}$.  
Using the notation of~\ref{para3}, the 
previous paragraph shows that the following diagram
\[ \xymatrix{
 \sdc{L^{\dagger}} & & \cdots & & \sdc{K_1^{\dagger}} & & \sdc{K^{\dagger}} & \\
\sdc{L^{\dagger}} \ar[u]_= 
& \sdc{L_{n-1}^{\dagger}} \ar[ul]\ar[ur] & & 
\sdc{L_1^{\dagger}} \ar[ul]\ar[ur] & & \sdc{L_0^{\dagger}} \ar[ul]\ar[ur] & 
\sdc{K^{\dagger}} \ar[u]_=
} \]
is a route $\route^{\dagger}$ from $\sdc{L^{\dagger}}$ 
to $\sdc{K^{\dagger}}$ with 
$\len_R(\route^{\dagger})=\len_R(\route)$.  This explains~\eqref{eqn2} 
and completes the proof.
\end{proof}

The next result is Theorem C from the introduction.  It will yield an
answer to the following:  In Proposition~\ref{prop4}, what is implied by the existence of
a fixed point for $\Delta$?  See Corollary~\ref{fixed02} for the answer.

\begin{thm} \label{fixed01}
Let $R$ be a local ring and $K,L$ homologically 
bounded below and degreewise finite 
$R$-complexes. 
If $K\lotimes_RK\lotimes_R L$ is 
semidualizing, then $K\sim R$.
\end{thm}

\begin{proof}
To keep the bookkeeping simple,
apply appropriate suspensions to $K$ and $L$ to assume 
$\inf(K)=0=\inf(L)$.  
Let $P$ and $Q$ be minimal projective resolutions of $K$ and
$L$, respectively; in particular, $P_0,Q_0\neq 0$.  
Then 
$K\lotimes_R K\lotimes_R L
\simeq P\otimes_R P\otimes_R Q$
is semidualizing, and so the 
homothety morphism is a 
quasi-isomorphism.
\[ R\xra{\simeq}
\Hom_R(P\otimes_R P\otimes_R Q,P\otimes_R P\otimes_R Q) \]
Here is the crucial point.  For complexes $X,Y$, let $\theta_{XY}\colon 
X\otimes_R Y\to Y\otimes_R X$ be the natural isomorphism.  This gives a 
cycle 
\[ \theta_{PP}\otimes_R Q\in
\Hom_R(P\otimes_R P\otimes_R Q,P\otimes_R P\otimes_R Q) \]
and therefore, 
there exists $u\in R$ such that the homothety 
$\mu_u\colon P\otimes_R P\otimes_R Q\to 
P\otimes_R P\otimes_R Q$ is homotopic to $\theta_{PP}\otimes_R Q$.  

Set $\ol{P}=P\otimes_R k$ and 
$\ol{Q}=Q\otimes_R k$. 
The fact that 
$\theta_{PP}\otimes_R Q$ and $\mu_u$ are homotopic implies that the 
following morphisms are also homotopic.
\[ (\theta_{PP}\otimes_R Q)\otimes_R k,\mu_u\otimes_R k \colon
(P\otimes_RP\otimes_R Q)\otimes_R k\to (P\otimes_RP\otimes_R Q)\otimes_R k \]
Using the 
isomorphism $(P\otimes_RP\otimes_R Q)\otimes_R k\cong 
\ol{P}\otimes_k\ol{P}\otimes_k\ol{Q}$, we then 
deduce that the $k$-morphisms
\[ \theta_{\ol{P}\,\ol{P}}\otimes_k\ol{Q},\mu_{\ol{u}}\colon 
\ol{P}\otimes_k\ol{P}\otimes_k\ol{Q}
\to\ol{P}\otimes_k\ol{P}\otimes_k\ol{Q}\]
are homotopic as well.  The differential on $\ol{P}\otimes_k\ol{P}\otimes_k\ol{Q}$ 
is 
zero by the minimality of $P$ and $Q$, and it follows that 
$\theta_{\ol{P}\,\ol{P}}\otimes_k\ol{Q}$ and $\mu_{\ol{u}}$ are equal.

Set $n=\rank_k\ol{P}_0$.
We claim that $n=1$.  
Suppose that
$n>1$, and   
fix bases $x_1,\ldots,x_n\in \ol{P}_0$ 
and $y_1,\ldots,y_p\in \ol{Q}_0$.  The set 
\[ \{x_i\otimes x_j\otimes y_l\mid \text{$i,j\in\{1,\ldots,n\}$ and
$l\in\{1,\ldots,p\}$}\}\]
is a basis for $\ol{P}_0\otimes_k \ol{P}_0\otimes_k \ol{Q}_0$.  
However, the equality 
\[ x_1\otimes x_2\otimes y_1 - ux_2\otimes x_1\otimes y_1=0 \]
contradicts the linear independence.
Thus, $n\leq 1$ and since $\ol{P}_0\neq 0$ we have $n=1$.

Next, we show that $\ol{P}_i=0$ for $i>0$.  The equality 
$\theta_{\ol{P}\,\ol{P}}\otimes_k\ol{Q}=\mu_{\ol{u}}$
implies 
\[ 0=x\otimes x'\otimes y\pm ux'\otimes x\otimes y 
\in(\ol{P}_0\otimes_k\ol{P}_i\otimes_k\ol{Q}_0)\oplus
(\ol{P}_i\otimes_k\ol{P}_0\otimes_k\ol{Q}_0) \]
for each $x\in \ol{P}_0$ and $x'\in\ol{P}_i$ and $y\in \ol{Q}_0$.
Since $\ol{P}_0\neq 0$ and 
$$(\ol{P}_0\otimes_k\ol{P}_i\otimes_k\ol{Q}_0)\cap
(\ol{P}_i\otimes_k\ol{P}_0\otimes_k\ol{Q}_0)=0$$
this is impossible unless $\ol{P}_i=0$.

One concludes
that there is an isomorphism
$P\simeq R$, completing the proof.
\end{proof}

\begin{cor} \label{corfixed}
Let $R$ be a local ring and $\sdc{K},\sdc{L}\in\s(R)$.
\begin{enumerate}[\quad\rm(a)]
\item \label{item201}
If $\sdc{K}\tri\sdc{L}$ and $\sdc{L}\tri\sdc{\rhom_R(L,K)}$,
then $\sdc{L}=\sdc{K}$.
\item \label{item202}
If $\sdc{K}\tri\sdc{L}$ and $\sdc{\rhom_R(L,K)}\tri\sdc{L}$,
then $\sdc{L}=\sdc{R}$.
\end{enumerate}
In particular, if $\sdc{K}\trineq\sdc{L}\trineq\sdc{R}$, then $\sdc{L}$ and 
$\sdc{\rhom_R(L,K)}$ are not comparable in the ordering on $\s(R)$.
\end{cor}

\begin{proof}
\eqref{item201} If $\rhom_R(L,K)$ is $L$-reflexive and
$L$ is $K$-reflexive, then~\ref{gerko1} provides 
$$K\simeq \rhom_R(L,K)\lotimes_R\rhom_R(\rhom_R(L,K),L)\lotimes_R 
\rhom_R(L,K).$$
Theorem~\ref{fixed01} then yields  $\rhom_R(L,K)\sim 
R$ and thus the second isomorphism in the next sequence
while the first follows since $\sdc{K}\tri\sdc{L}$ and the third 
is standard.
$$L\simeq\rhom_R(\rhom_R(L,K),K)\sim\rhom_R(R,K)\simeq K$$

\eqref{item202} If $L$ is $\rhom_R(L,K)$-reflexive and
$K$-reflexive, then the isomorphism
$L\simeq\rhom_R(\rhom_R(L,K),K)$ 
with part~\eqref{item201} implies that $\rhom_R(L,K)\sim K$.  The 
desired isomorphism then follows from an application of $\rhom_R(-,K)$.

The final statement follows directly from parts~\eqref{item201} 
and~\eqref{item202}.
\end{proof}

In view of condition~\eqref{fixed000} of the next result
we note the following open 
question:  If $R$ is a local ring, must $\s(R)$ be a finite set?  The 
answer is known in very few cases.  
See~\cite{LWC:new} and~\cite{sather:divisor} for discussion of this question.

\begin{cor} \label{fixed02}
For a local ring $R$, the following conditions are equivalent:
\begin{enumerate}[\quad\rm(i)]
\item \label{fixed00}
$R$ is Gorenstein.
\item \label{fixed0}
$R$ admits a dualizing complex $D$ and a semidualizing 
complex $L$ such that $\sdc{\rhom_R(L,D)}=\sdc{L}$.
\item \label{fixed000}
$R$ admits a dualizing complex $D$, and $\s(R)$ is 
finite with odd cardinality.
\end{enumerate}
\end{cor}

\begin{proof}
The implications \eqref{fixed00}$\implies$\eqref{fixed0} and
\eqref{fixed00}$\implies$\eqref{fixed000} are
clear as, when $R$ is Gorenstein, $R$ is dualizing and 
$\s(R)=\{\sdc{R}\}$ by~\cite[(8.6)]{christensen:scatac}.  

\eqref{fixed0}$\implies$\eqref{fixed00}.
Let $D,L$ be as 
in~\eqref{fixed0}.
Corollary~\ref{corfixed}\eqref{item201} with $K=D$ provides the first 
and third isomorphisms in
the next sequence 
$$D\sim L\sim \rhom_R(L,D)\sim \rhom_R(D,D)\simeq R$$
while the others follow by hypothesis.  Thus,
$R$ 
is Gorenstein.

\eqref{fixed000}$\implies$\eqref{fixed0}.
Assume that $R$ admits a dualizing complex $D$ and that
$\s(R)$ is finite.  
If $\sdc{\rhom_R(L,D)}\neq\sdc{L}$ for all $\sdc{L}\in\s(R)$,
then 
$\s(R)$ is the 
disjoint 
union of subsets of the form $\{\sdc{L},\sdc{\rhom_R(L,D)}\}$, each of 
which has two distinct elements.  Thus, $\s(R)$ has even cardinality, 
completing the proof.
\end{proof}

Here is Theorem B from the introduction.
For a specific construction of 
noncomparable semidualizing complexes, see Corollary~\ref{noncom}.
For $\sdc{K}\in\s(R)$ and $\delta>0$, set
$B(\sdc{K},\delta)=\{\sdc{L}\in\s(R)\mid\dist(\sdc{K},\sdc{L})<\delta\}$.

\begin{thm} \label{prop201}
For a local ring $R$ the following conditions are equivalent:
\begin{enumerate}[\quad\rm(i)]
\item \label{item207}
There exist elements of $\s(R)$ that are not comparable.
\item \label{item208}
$\s(R)$ has cardinality at least 3.
\item \label{item209}
There exists $\sdc{K}\in\s(R)$ and $\delta>0$ such that the open 
ball $B(\sdc{K},\delta)$ satisfies $\{\sdc{K}\}\subsetneq 
B(\sdc{K},\delta)\subsetneq\s(R)$.
\end{enumerate}
\end{thm}

\begin{proof}
\eqref{item208}$\implies$\eqref{item207}
Fix distinct elements 
$\sdc{K}, \sdc{L}, 
\sdc{M}\in\s(R)$.  Without loss of generality, assume that $\sdc{M}=\sdc{R}$.
Suppose that every two 
elements in $\s(R)$ are comparable.  
The elements $\sdc{K}, \sdc{L}, 
\sdc{R}$ can be reordered to assume that $\sdc{K}\trineq \sdc{L}\trineq 
\sdc{R}$, and we are done by 
Corollary~\ref{corfixed}.

\eqref{item209}$\implies$\eqref{item208} Let $\sdc{K}\in\s(R)$ and 
$\delta>0$ be such that $\{\sdc{K}\}\subsetneq 
B(\sdc{K},\delta)\subsetneq\s(R)$.  Fixing $\sdc{L}\in\s(R) 
\smallsetminus B(\sdc{K},\delta)$ and $\sdc{M}\in 
B(\sdc{K},\delta) \smallsetminus \{\sdc{K}\}$ provides at least three 
distinct elements of $\s(R)$: $\sdc{K},\sdc{L},\sdc{M}$.

\eqref{item207}$\implies$\eqref{item209} Fix two noncomparable elements
$\sdc{K},\sdc{L}\in\s(R)$, and let $\route$ be a route in 
$\Graph(R)$ from $\sdc{K}$ to $\sdc{L}$ such that 
$\len(\route)<\dist(\sdc{K},\sdc{L})+\frac{1}{2}$.
Using the notation of~\ref{para3} for $\route$, there exists an 
integer $i$ between $0$ and $n$ such that either 
$\sdc{K}\neq\sdc{K_i}$ or 
$\sdc{K}=\sdc{K_i}\neq\sdc{L_i}$, and we let $i_0$ 
denote the smallest 
such integer.
If $\sdc{K}\neq\sdc{K_{i_0}}$, then 
$\sdc{K_{i_0}}\trineq\sdc{K}$.  In this event, 
$\route$  can be factored as the concatenation
$\route_1\route_2$ as in the following 
diagram.
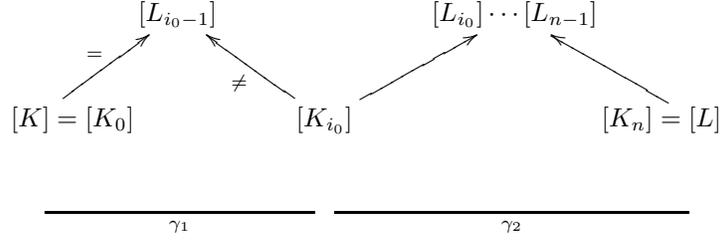
\[ \xymatrix{
 & \sdc{L_{i_0-1}} & & \sdc{L_{i_0}} \cdots  \sdc{L_{n-1}} & \\
\sdc{K}=  \sdc{K_0}\hspace{-1cm} \ar[ur]^= & & 
\sdc{K_{i_0}}  \ar[ul]^{\neq}\ar[ur] & & 
\hspace{-1cm}\sdc{K_n}=  \sdc{L} \ar[ul] \\
\vspace{-5mm}
\, & & \, \ar@{-}[ll]^{\route_1}\ar@{-}[rr]_{\route_2} & & \,
} \]
Since $\sdc{K},\sdc{L}$ are not comparable, it follows that 
$\sdc{K_{i_0}}\neq\sdc{L}$ and so $\len(\route_2)>\frac{1}{2}$ 
by Lemma~\ref{para1}.  With Theorem~\ref{prop2} this provides 
(1) in the following sequence
$$\dist(\sdc{K},\sdc{K_{i_0}})+\textstyle\frac{1}{2}\displaystyle
\stackrel{(1)}{<}\len(\route_1)+\len(\route_2)
\stackrel{(2)}{=}\len(\route)
\stackrel{(3)}{<}\dist(\sdc{K},\sdc{L})+\textstyle\frac{1}{2}$$
while (2) is by~\ref{subprop3} and (3) is from the choice 
of $\route$.
In particular, $\dist(\sdc{K},\sdc{K_{i_0}})<\dist(\sdc{K},\sdc{L})$.
Fixing $\delta$ such that $\dist_R(\sdc{K},\sdc{K_{i_0}})<\delta<\dist_R(\sdc{K},\sdc{L})$, 
one has 
$\sdc{K_{i_0}}\in B(\sdc{K},\delta)\smallsetminus\{\sdc{K}\}$ and
$\sdc{L}\in\s(R)\smallsetminus B(\sdc{K},\delta)$, giving the 
desired proper containments.

If $\sdc{K}=\sdc{K_{i_0}}\neq\sdc{L_{i_0}}$, then similar 
reasoning shows that, with
a choice of $\delta$ such that $\dist_R(\sdc{K},\sdc{L_{i_0}})<\delta<\dist_R(\sdc{K},\sdc{L})$, 
one has 
$\sdc{L_{i_0}}\in B(\sdc{K},\delta)\smallsetminus\{\sdc{K}\}$ and
$\sdc{L}\in\s(R)\smallsetminus B(\sdc{K},\delta)$.
\end{proof}

\section{Behavior of the metric under change of rings} \label{sec3}

\noindent \emph{In this section, let $\vf\colon R\to S$ be a local homomorphism of finite
flat dimension.}

\begin{para} \label{bc}
\textbf{Base change:}
The homomorphism $\vf$ induces a well-defined injective map
$$\s_{\vf}\colon\s(R)\to \s(S) \qquad\text{given by}\qquad \sdc{K}\mapsto\sdc{K\lotimes_R S}$$
by~\cite[(4.5),(4.9)]{frankild:appx}, and $\sdc{K}\tri\sdc{L}$ if and only if
$\sdc{K\lotimes_R S}\tri\sdc{L\lotimes_R S}$ by~\cite[(4.8)]{frankild:appx}.
When $\sdc{K}\tri\sdc{L}$, one has
$P^S_{\rhom_S(L\lotimes_R S,K\lotimes_R S)}(t)
=P^R_{\rhom_R(L,K)}(t)$
from~\cite[(6.15)]{frankild:appx},
providing the second equality in the following sequence
$$\dist_{S}(\sdc{K\lotimes_R S},\sdc{L\lotimes_R S})
=\prox_{S}(\sdc{K\lotimes_R S},\sdc{L\lotimes_R S}) 
=\prox_R(\sdc{K},\sdc{L})  
=\dist_R(\sdc{K},\sdc{L}) $$
while the other equalities are from Theorem~\ref{prop2}.
\end{para}

Next we show that the metric is nonincreasing under $\s_{\vf}$;  
we do not know of an example where it decreases.
For instance, equality holds in each case  where $\s(S)$
is completely determined in~\cite{sather:divisor}.
When $\s_{\vf}$ is surjective, the result states that $\s_{\vf}$ is an isometry.

\begin{prop}  \label{prop5} 
Let $\vf\colon R\to S$ be a local homomorphism of finite flat
dimension.  For $\sdc{K},\sdc{L}\in\s(R)$ there is an inequality
\[ \dist_{S}(\sdc{K\lotimes_R S},\sdc{L\lotimes_R S})\leq 
\dist_R(\sdc{K},\sdc{L}) \]
with equality when $\sdc{K}\tri\sdc{L}$ or when 
$\s_{\vf}$ is surjective, e.g., if $R$ is complete and $\vf$ is surjective with
kernel generated by an $R$-sequence.
\label{thm10}
\end{prop}

\begin{proof}
Fix $\sdc{K},\sdc{L}$ in $\s(R)$.
When $\sdc{K}\tri\sdc{L}$, the equality  is in~\ref{bc}.
In general, let $\route$ 
be a route from $\sdc{K}$ to $\sdc{L}$ in $\Graph(R)$.
Using the notation of~\ref{para3}, the 
diagram
\[ \xymatrix{
 & \sdc{L_0\lotimes_R S} &  & \sdc{L_{n-1}\lotimes_R S} \\
\sdc{K_0\lotimes_R S} \ar[ur] & & \cdots \ar[ul]\ar[ur]
& & \sdc{K_n\lotimes_R S} \ar[ul]
} \]
is a route $\route\lotimes_R S$ from 
$\sdc{K\lotimes_R S}$ to $\sdc{L\lotimes_R S}$
in $\Graph(S)$ by~\ref{bc}  and 
\[ \dist_{S}(\sdc{K\lotimes_R S},\sdc{L\lotimes_R S}) 
\leq \len_S(\route\lotimes_R S) =\len_R(\route). \]
Since this is true for every route $\route$, the desired inequality 
now follows.

When $\s_{\vf}$ is surjective, 
the above analysis along with~\ref{bc} 
implies that the routes from $\sdc{K}$ to $\sdc{L}$ are in length-preserving 
bijection with those from 
$\sdc{K\lotimes_R S}$ to $\sdc{L\lotimes_R S}$ and so
$\dist_{S}(\sdc{K\lotimes_R S},\sdc{L\lotimes_R S})= 
\dist_R(\sdc{K},\sdc{L})$.
When $R$ is complete and $\vf$ is surjective with
kernel generated by an $R$-sequence,
the surjectivity of $\s_{\vf}$  follows from~\cite[(4.5)]{frankild:appx}
and~\cite[(3.2)]{yoshino}.  
\end{proof}

Using~\cite[(2.5),(3.16)]{christensen:scatac} 
and the inequality $P^{R_{\p}}_{X_{\p}}(t)\preccurlyeq 
P^R_X(t)$ with~\ref{product}\eqref{item01},
the proof of the previous result
yields the following.  Example~\ref{decreasing}  shows that inequality
may be strict or not.

\begin{prop} \label{thm8}
For $\p\in\spec(R)$ and $\sdc{K},\sdc{L}\in\s(R)$, there is an inequality
\[ 
\phantom{\qquad\,\quad\qquad\qquad\qquad}
\dist_{R_{\p}}(\sdc{K_{\p}},\sdc{L_{\p}})\leq\dist_R(\sdc{K},\sdc{L}). 
\phantom{\qquad\,\quad\qquad\qquad\qquad}\qed 
\]
\end{prop}

\begin{para} \label{cc}
\textbf{Cobase change:}
A \emph{Gorenstein factorization} of $\vf$ is a diagram of local homomorphisms
$R\xra{\dot\vf}R'\xra{\vf'}S$
such that $\vf=\vf'\dot\vf$, $\vf'$ is surjective, and $\dot{\vf}$ is flat with 
Gorenstein closed fibre.  
Homomorphisms admitting Gorenstein factorizations 
exist in profusion, e.g., if $\vf$ is essentially of finite type or if
$S$ is complete;
see~\cite[(1.1)]{avramov:solh}. 

Assume that $\vf$ admits a Gorenstein factorization 
as above and set $d=\depth(\dot{\vf})$.  
For each homologically finite
complex
of $R$-modules $X$, 
set 
\[ \cbc{X}{\vf}=\shift^{d}\rhom_{R'}(S,X\lotimes_R R'). \]
It is shown in~\cite[(6.5),(6.12)]{frankild:appx} that this
is independent of Gorenstein factorization and that 
the following assignment is well-defined
and injective.
$$\s^{\vf}\colon\s(R)\to\s(S) \qquad\text{given by}\qquad
\sdc{K}\mapsto\sdc{\cbc{K}{\vf}}$$
One has $\sdc{K}\tri\sdc{L}$ if and only if
$\sdc{\cbc{K}{\vf}}\tri\sdc{\cbc{L}{\vf}}$ by~\cite[(6.11)]{frankild:appx}.
When $\sdc{K}\tri\sdc{L}$, one has
$P^S_{\rhom_S(\cbc{L}{\vf},\cbc{K}{\vf})}(t)
=P^R_{\rhom_R(L,K)}(t)$
from~\cite[(6.15)]{frankild:appx},
providing the second equality in the following sequence
while the others are from Theorem~\ref{prop2}.
$$\dist_{S}(\sdc{\cbc{K}{\vf}},\sdc{\cbc{L}{\vf}})
=\prox_{S}(\sdc{\cbc{K}{\vf}},\sdc{\cbc{L}{\vf}}) 
=\prox_R(\sdc{K},\sdc{L})  
=\dist_R(\sdc{K},\sdc{L}) $$

\label{dipsy}
One has $\sdc{K}\tri \sdc{L}$ if and only if 
$\sdc{\cbc{K}{\vf}}\tri\sdc{L\lotimes_RS}$ by~\cite[(6.13)]{frankild:appx}.
When  $\sdc{K}\tri \sdc{L}$ one has
$P^S_{\rhom_S(L\lotimes_RS,\cbc{K}{\vf})}(t)=
P^R_{\rhom_R(L,K)}(t)I_{\vf}(t)$ from~\cite[(6.15)]{frankild:appx}, and
using~\ref{product}\eqref{item02}
and Theorem~\ref{prop2} as above there is an equality
\begin{equation}\label{dippy0}\tag{$\dagger$}
\dist_{S}(\sdc{\cbc{K}{\vf}},\sdc{L\lotimes_R S})=
\max\{\dist_R(\sdc{K},\sdc{L}),\icurv(\vf)\}.
\end{equation}
Example~\ref{nontrivial} shows that this formula can fail if
$\sdc{K}\not\tri \sdc{L}$.  See Corollary~\ref{cor4} for the general case.
The  case $\sdc{K}=\sdc{L}$ yields 
\begin{equation}\label{dippy}\tag{$\ddagger$}
\dist_{S}(\sdc{\cbc{K}{\vf}},\sdc{K\lotimes_R S})=
\icurv(\vf)
\end{equation}
and thus $\vf$ is Gorenstein at $\n$ if and only if 
$\s_{\vf}=\s^{\vf}$; see~\ref{paraRadii}. 
\end{para}

As with Proposition~\ref{prop5} we do not know if the 
next inequality can be strict.  

\begin{prop} \label{prop6}
Let $\vf\colon R\to S$ be a local homomorphism of finite flat
dimension admitting a 
Gorenstein factorization. 
For $\sdc{K},\sdc{L}$ in $\s(R)$,
there is an inequality
$$\dist_{S}(\sdc{\cbc{K}{\vf}},
  \sdc{\cbc{L}{\vf}})\leq\dist_R(\sdc{K},\sdc{L}) $$
with equality when $\sdc{K}\tri\sdc{L}$.  Equality also holds when $\s^{\vf}$ is surjective,
in which case $\vf$ is Gorenstein at $\n$.
\end{prop}

\begin{proof}
When $\sdc{K}\tri \sdc{L}$, 
the equality is in~\ref{cc}.
As in the proof of Proposition~\ref{prop5}, 
for arbitrary $\sdc{K},\sdc{L}$, let $\route$ be a route from $\sdc{K}$ 
to $\sdc{L}$ in $\Graph(R)$, with
the notation of~\ref{para3}. The 
following diagram
is a route $\cbc{\route}{\vf}$ from $\sdc{\cbc{K}{\vf}}$ 
to $\sdc{\cbc{L}{\vf}}$ 
\[ \xymatrix{
 & \sdc{\cbc{L_0}{\vf}} &  & 
\sdc{\cbc{L_{n-1}}{\vf}} & \\
\sdc{\cbc{K_0}{\vf}} \ar[ur] & & 
\cdots  \ar[ul]\ar[ur] & & \sdc{\cbc{K_n}{\vf}} \ar[ul]
} \]
with $\dist_{S}(\sdc{\cbc{K}{\vf}},\sdc{\cbc{L}{\vf}}) \leq
\len_{S}(\cbc{\route}{\vf})=\len_R(\route)$.

If $\s^{\vf}$ is surjective, then there exists $\sdc{K}\in\s(R)$ such that
$\sdc{\cbc{K}{\vf}}=\sdc{S}$.  With $L=R$ in~\ref{cc} equation~\eqref{dippy0},
one has $\icurv(\vf)=0$, so $\vf$ is Gorenstein at $\n$ and $\s_{\vf}=\s^{\vf}$.
The equality now follows from Proposition~\ref{prop5}.  
\end{proof}

By equation~\eqref{dippy0} of~\ref{dipsy} and Example~\ref{nontrivial},
the next inequality may be strict or not.

\begin{cor} \label{cor4}
Let $\vf\colon R\to S$ be a local homomorphism of finite flat
dimension admitting a 
Gorenstein factorization. 
For $\sdc{K},\sdc{L}$ in $\s(R)$  
there is an inequality
$$\dist_{S}(\sdc{\cbc{K}{\vf}},\sdc{L\lotimes_R S})\leq
\icurv(\vf)+\dist_R(\sdc{K},\sdc{L}).$$
\end{cor}

\begin{proof}
Use the triangle inequality,
\ref{dipsy} equation~\eqref{dippy}, and Proposition~\ref{prop5}
\begin{align*}
\dist_{S}(\sdc{\cbc{K}{\vf}},\sdc{L\lotimes_R S}) 
&\leq\dist_{S}(\sdc{\cbc{K}{\vf}},\sdc{K\lotimes_R S})
+\dist_{S}(\sdc{K\lotimes_R S},\sdc{L\lotimes_R S}) \\
&\leq   \icurv(\vf)+\dist_R(\sdc{K},\sdc{L}). \qedhere
\end{align*}
\end{proof}

See Example~\ref{nontrivial} for a special case of the
final statement of the next result.

\begin{cor} \label{noncom}
Let $\vf\colon R\to S$ be a 
local homomorphism of finite flat 
dimension admitting a Gorenstein factorization 
and fix $\sdc{K},\sdc{L}\in \s(R)$ with $\sdc{K}\tri\sdc{L}$.  
\begin{enumerate}[\quad\rm(a)]
\item \label{item204}
If $\sdc{L\lotimes_R S}\tri\sdc{\cbc{\rhom_R(L,K)}{\vf}}$, 
then $\vf$ is Gorenstein at 
$\n$ and $\sdc{K}=\sdc{L}$.
\item \label{item205}
If $\sdc{\cbc{\rhom_R(L,K)}{\vf}}\tri\sdc{L\lotimes_R S}$, then 
$\sdc{L}=\sdc{R}$.
\end{enumerate}
In particular, if $\sdc{K}\neq\sdc{R}$ and $\vf$ is not Gorenstein 
at $\n$, then 
the elements $\sdc{K\lotimes_R S}$ and
$\sdc{\cbc{R}{\vf}}$ are noncomparable in the 
ordering on $\s(S)$.
\end{cor}

\begin{proof}
The assumption $\sdc{K}\tri\sdc{L}$ implies $\sdc{\cbc{K}{\vf}}\tri\sdc{L\lotimes_R S}$
by~\ref{dipsy}, and~\cite[(6.9)]{frankild:appx} provides an
isomorphism
$\rhom_S(L\lotimes_R S,\cbc{K}{\vf})\simeq\cbc{\rhom_R(L,K)}{\vf}$.

If 
$\sdc{L\lotimes_R S}\tri\sdc{\cbc{\rhom_R(L,K)}{\vf}}
=\sdc{\rhom_S(L\lotimes_R S,\cbc{K}{\vf})}$, 
then 
Corollary~\ref{corfixed}\eqref{item201} 
implies $\sdc{L\lotimes_R S}=\sdc{\cbc{K}{\vf}}$.
Equation~\eqref{dippy0} in~\ref{dipsy} yields the conclusion for part~\eqref{item204}.
Part~\eqref{item205} follows similarly from 
Corollary~\ref{corfixed}\eqref{item202},  
and the final statement is a consequence of~\eqref{item204} and 
~\eqref{item205} using $K=L$.
\end{proof}

\section{Examples} \label{sec4}

This section consists of specific computations of distances in $\s(R)$.
We begin with a simple example upon which the others are built.  It 
shows, in particular, that although the diameter of the metric space 
$\s(R)$ is finite by Proposition~\ref{para5}, it can be arbitrarily 
large.  Here, the diameter of
$\s(R)$ is
\[ \operatorname{diam}(\s(R))=\sup\{\dist_R(\sdc{K},\sdc{L})\mid 
\sdc{K},\sdc{L}\in\s(R)\}. \]

\begin{ex} \label{square0}
Assume that $\m^2=0$.  In particular, $R$ is Cohen-Macaulay, so each 
semidualizing complex is, up to shift, isomorphic to a module 
by~\cite[(3.7)]{christensen:scatac}.  Since $R$ is 
Artinian, it admits a dualizing module $D$ by~\ref{homothety}.  
The set $\s(R)$ contains at most two distinct elements, namely 
$\sdc{R}$ and $\sdc{D}$:  If $K$ is a nonfree semidualizing module, 
then any syzygy module from a minimal free resolution
of $K$ is a nonzero $k$-vector space that is $K$-reflexive, implying 
that $K$ is dualizing by~\cite[(8.4)]{christensen:scatac}.

The elements $\sdc{R}$ and $\sdc{D}$ are distinct if and only if 
$R$ is non-Gorenstein.  When these conditions hold, the previous argument 
shows that $\curv_R(D)=\curv_R(k)$.  A 
straightforward computation of the minimal free resolution of $k$ 
shows that 
\[ P^R_k(t)=\sum_{n=0}^{\infty}r^nt^n=1/(1-rt) \]
where 
$r=\edim(R)=\rank_k(\m/\m^2)$.  In particular, 
\[ \dist_R(\sdc{R},\sdc{D})=\curv_R(D)=\curv_R(k)=r \]  
and thus 
$\operatorname{diam}(\s(R))=r$.  
The trivial extension $k\ltimes k^r$ gives an 
explicit example.
\end{ex}

We now give a particular example of 
the construction from Corollary~\ref{noncom}
which has the added
benefit of being an example where we can completely describe the 
structure of the metric space $\s(S)$.  Note that this process can be 
iterated.

\begin{ex} \label{nontrivial}
Fix integers $r,s\geq 2$ and a field $k$.  Let 
$R=k\ltimes k^r$ and $S=R\ltimes R^s$.  The natural map $\vf\colon 
R\to S$ is flat and local with closed fibre $\ol{S}\cong k\ltimes k^s$.  
Since $R$ is Artinian it admits a dualizing module  $D$ by~\ref{homothety}.  
By~\cite[(4.7)]{sather:divisor}, the set $\s(S)$ consists of the 
four distinct elements $\sdc{S},\sdc{D\lotimes_R S}, 
\sdc{\cbc{R}{\vf}},\sdc{\cbc{D}{\vf}}$. 
The next Poincar\'{e} series and 
curvatures are computed using~\ref{product}\eqref{item02}, Example~\ref{square0}, 
and~\cite[(6.10),(6.15)]{frankild:appx}. 
\begin{align*}
&P^S_{D\lotimes_R S}
=I^R_R(t) 
&  & \curv_S(D\lotimes_R S) = r \\
&P^S_{\cbc{R}{\vf}}(t)
= I^{\ol{S}}_{\ol{S}}(t) 
&  & \curv_S(\cbc{R}{\vf}) = s \\
&P^S_{\cbc{D}{\vf}}(t)
=I^R_R(t)I^{\ol{S}}_{\ol{S}}(t)
&  & \curv_R(\cbc{D}{\vf})  =\max\{ r,s\} \\
&P^S_{\rhom_S(D\lotimes_R S,\cbc{D}{\vf})}(t)
=I^{\ol{S}}_{\ol{S}}(t) 
&  & \curv_S(\rhom_S(D\lotimes_R S,\cbc{D}{\vf})) = s \\
&P^S_{\rhom_S(\cbc{R}{\vf},\cbc{D}{\vf})}(t)
=I^R_R(t) 
& & \curv_S(\rhom_S(\cbc{R}{\vf},\cbc{D}{\vf})) = r
\intertext{With Theorem~\ref{prop2},
this gives 
the following distance computations.}
&\dist_S(\sdc{S},\sdc{D\lotimes_R S}) =r 
&&\dist_S(\sdc{S},\sdc{\cbc{R}{\vf}}) =s \\
&\dist_R(\sdc{S},\sdc{\cbc{D}{\vf}}) =\max\{ r,s\} 
&&\dist_S(\sdc{D\lotimes_R S},\sdc{\cbc{D}{\vf}}) =s \\
&\dist_S(\sdc{\cbc{R}{\vf}},\sdc{\cbc{D}{\vf}}) =r
\end{align*}
This provides the lengths of the edges in the following sketch of $\Graph(S)$
\[ \xymatrix{ 
 & \sdc{S} & \\
\sdc{D\lotimes_R S} \ar[ur]^r & & & \sdc{\cbc{R}{\vf}} \ar[ull]^s \\
& & \sdc{\cbc{D}{\vf}} \ar[ull]^s \ar[ur]^r \ar[uul]^{\max\{ r,s\}}
} \]
while~\cite[(4.7)]{sather:divisor} implies that this is a complete 
description of $\Graph(S)$.  Thus, the remaining distance
is computed readily:
$$\dist_S(\sdc{\cbc{R}{\vf}},\sdc{D\lotimes_R S}) =r+s.$$
In particular, the open ball in $\s(S)$
of radius $r+1$ centered at $\sdc{\cbc{R}{\vf}}$ contains 
$\sdc{\cbc{D}{\vf}}\neq \sdc{\cbc{R}{\vf}}$ and does not contain  
$\sdc{D\lotimes_R S}$.   
Furthermore, this shows that equality can hold in Corollary~\ref{cor4}.
\end{ex}

Finally, we show that the metric may or may not decrease 
after localizing.

\begin{ex} \label{decreasing}
Let 
$R$ be a non-Gorenstein ring 
with dualizing complex $D$ and $\p$ a prime 
ideal such that $R_{\p}$ is Gorenstein, e.g.,
$R=k[\![X,Y]\!]/(X^2,XY)$ and $\p=(X)R$.
Then $D_{\p}\sim R_{\p}$, implying
\[
\dist_{R_{\p}}(\sdc{R_{\p}},\sdc{D_{\p}})=0 
<\dist_R(\sdc{R},\sdc{D}). \]
On the other hand, let $S=k[\![X,Y,Z]\!]/(X,Y)^2$ with dualizing module 
$E$, and $\q=(X,Y)S$;  then the computations in
Example~\ref{square0} give
\[ \dist_{S_{\q}}(\sdc{S_{\q}},\sdc{E_{\q}}) 
=\curv_{S_{\q}}(E_{\q})=2 \]
while Proposition~\ref{thm10} yields the first equality in the following 
sequence
\begin{align*}
\dist_S(\sdc{S},\sdc{E})
&=\dist_{S/(Z)}(\sdc{S\lotimes_S S/(Z)},\sdc{E\lotimes_S S/(Z)}) \\
&=\curv_{S/(Z)}(E\lotimes_S S/(Z))=2 
\end{align*}
and the last equality follows from Example~\ref{square0}
since $S/(Z)\cong k[X,Y]/(X,Y)^2$.
\end{ex}

\section*{Acknowledgments}
 
A.F.~is grateful to the Department of Mathematics at the 
University of Illinois at Urbana-Champaign for its hospitality while 
much of this research was conducted.  S.S.-W.~is similarly 
grateful to the Institute for Mathematical 
Sciences at the University of Copenhagen.
Both authors thank Luchezar Avramov,
Lars Winther Christensen, E.~Graham Evans Jr.,
Hans-Bj\o rn Foxby, Alexander Gerko, Phillip Griffith,
Henrik Holm, Srikanth Iyengar, and Paul Roberts for 
stimulating conversations and helpful comments about this 
research, and the anonymous referee for improving the presentation.  

\providecommand{\bysame}{\leavevmode\hbox to3em{\hrulefill}\thinspace}
\providecommand{\MR}{\relax\ifhmode\unskip\space\fi MR }
\providecommand{\MRhref}[2]{
  \href{http://www.ams.org/mathscinet-getitem?mr=#1}{#2}
}
\providecommand{\href}[2]{#2}


\begin{thebibliography}{10}

\bibitem{auslander:adgeteac}
M.~Auslander, \emph{Anneaux de {G}orenstein, et torsion en alg\`ebre
  commutative}, S\'eminaire d'Alg\`ebre Commutative dirig\'e par Pierre Samuel,
  vol. 1966/67, Secr\'etariat math\'ematique, Paris, 1967. \MR{37 \#1435}

\bibitem{auslander:smt}
M.~Auslander and M.~Bridger, \emph{Stable module theory}, Memoirs of the
  American Mathematical Society, No. 94, American Mathematical Society,
  Providence, R.I., 1969. \MR{42 \#4580}

\bibitem{avramov:mwer}
L.~L. Avramov, \emph{Modules with extremal resolutions}, Math. Res. Lett.
  \textbf{3} (1996), 319--328. \MR{97f:13020}

\bibitem{avramov:ifr}
\bysame, \emph{Infinite free resolutions}, Six lectures on commutative algebra
  (Bellaterra, 1996), Progr. Math., vol. 166, Birkh\"auser, Basel, 1998,
  pp.~1--118. \MR{99m:13022}

\bibitem{avramov:rhafgd}
L.~L. Avramov and H.-B. Foxby, \emph{Ring homomorphisms and finite {G}orenstein
  dimension}, Proc. London Math. Soc. (3) \textbf{75} (1997), no.~2, 241--270.
  \MR{98d:13014}

\bibitem{avramov:solh}
L.~L. Avramov, H.-B. Foxby, and B.~Herzog, \emph{Structure of local
  homomorphisms}, J. Algebra \textbf{164} (1994), 124--145. \MR{95f:13029}

\bibitem{avramov:bsolrhoffd}
L.~L. Avramov, H.-B.Foxby, and J.~Lescot, \emph{Bass series of local ring
  homomorphisms of finite flat dimension}, Trans. Amer. Math. Soc. \textbf{335}
  (1993), no.~2, 497--523. \MR{93d:13026}

\bibitem{christensen:gd}
L.~W. Christensen, \emph{Gorenstein dimensions}, Lecture Notes in Mathematics,
  vol. 1747, Springer-Verlag, Berlin, 2000. \MR{2002e:13032}

\bibitem{christensen:scatac}
\bysame, \emph{Semi-dualizing complexes and their {A}uslander categories},
  Trans. Amer. Math. Soc. \textbf{353} (2001), no.~5, 1839--1883.
  \MR{2002a:13017}

\bibitem{LWC:new}
L.~W.~Christensen and S.~Sather-Wagstaff, \emph{Descent of semidualizing complexes 
for rings with the approximation property}, in preparation.

\bibitem{foxby:hacr}
H.-B. Foxby, \emph{Hyperhomological algebra \& commutative rings}, in
  preparation.

\bibitem{foxby:gmarm}
\bysame, \emph{Gorenstein modules and related modules}, Math. Scand.
  \textbf{31} (1972), 267--284 (1973). \MR{48 \#6094}

\bibitem{frankild:appx}
A.~Frankild and S.~Sather-Wagstaff, \emph{Reflexivity and ring homomorphisms of
  finite flat dimension}, Comm.~Algebra, to appear,
  \texttt{arXiv:math.AC/0508062}.

\bibitem{gelfand:moha}
S.~I. Gelfand and Y.~I. Manin, \emph{Methods of homological algebra},
  Springer-Verlag, Berlin, 1996. \MR{2003m:18001}

\bibitem{gerko:sdc}
A.~Gerko, \emph{On the structure of the set of semidualizing complexes},
  Illinois J. Math. \textbf{48} (2004), no.~3, 965--976. \MR{2114263}

\bibitem{golod:gdagpi}
E.~S. Golod, \emph{{$G$}-dimension and generalized perfect ideals}, Trudy Mat.
  Inst. Steklov. \textbf{165} (1984), 62--66, Algebraic geometry and its
  applications. \MR{85m:13011}

\bibitem{hartshorne:rad}
R.~Hartshorne, \emph{Residues and duality}, Lecture Notes in Mathematics, No.
  20, Springer-Verlag, Berlin, 1966. \MR{36 \#5145}

\bibitem{sather:divisor}
S.~Sather-Wagstaff, \emph{Semidualizing modules and the divisor class group},
  preprint \texttt{arXiv:math.AC/0404361}.

\bibitem{verdier:cd}
J.-L. Verdier, \emph{Cat\'{e}gories d\'{e}riv\'{e}es}, SGA 4$\frac{1}{2}$,
  Springer-Verlag, Berlin, 1977, Lecture Notes in Mathematics, Vol. 569,
  pp.~262--311. \MR{57 \#3132}

\bibitem{verdier:1}
\bysame, \emph{Des cat\'egories d\'eriv\'ees des cat\'egories ab\'eliennes},
  Ast\'erisque (1996), no.~239, xii+253 pp. (1997), With a preface by Luc
  Illusie, Edited and with a note by Georges Maltsiniotis. \MR{98c:18007}

\bibitem{yoshino}
Y.~Yoshino, \emph{The theory of {L}-complexes and weak liftings of complexes},
  J. Algebra \textbf{188} (1997), no.~1, 144--183. \MR{98i:13024}

\end{thebibliography}
\end{document}